\documentclass[a4paper]{article}
\usepackage{latexsym,amsfonts,amscd,amssymb,enumerate,amstext,amsmath,dsfont,amsthm}

\usepackage[all]{xy}

\newtheorem{thm}{Theorem}[section]
\newtheorem{prop}[thm]{Proposition}
\newtheorem{lem}[thm]{Lemma}
\newtheorem{cor}[thm]{Corollary}

\newtheorem{defn}[thm]{Definition}
\newtheorem{exmp}[thm]{Example}

\newcommand{\cata}{\mathcal{A}}
\newcommand{\catb}{\mathcal{B}}
\usepackage{setspace}
\usepackage{geometry}
\title{Adjoint functors and  triangulated categories}

\author{M. Grime}
\begin{document}
\doublespacing
\maketitle

\begin{abstract}
We give a construction of triangulated categories as quotients of exact categories where the subclass of objects sent to zero is defined by a triple of functors. This includes the cases of homotopy and stable module categories. These categories naturally fit into a framework of relative derived categories, and once we prove that there are decent resolutions of complexes, we are able to prove many familiar results in homological algebra.
\end{abstract}

\section{Introduction}

One of the zeniths of modular representation theory is Green's theory of vertices and sources, which gives, when properly formulated, a result about the existence of equivalences of additive categories. This relies on the notion that if one picks a subgroup (or collection of subgroups) of a finite group $G$ then there is a well-understood notion of projectivity relative to that subgroup. We may form the quotient category by quotienting out by these relatively projective objects and this can easily be shown to yield a triangulated category. More recently Okuyama in his unpublished manuscript \cite{okuyama} generalized the notion of relative projectivity to the so-called `with respect to a module' variety. The corresponding quotient is triangulated, as can be seen by appeal to Happel's theorem in \cite{happel}, and is usually called the $W$-stable category where $W$ is the module in question. 

Carlson, Peng and Wheeler, in \cite{virtual} further define a generalization of the transfer map in cohomology and use this to define a class of modules they call virtually $W$-projective. These objects form a thick subcategory of the $W$-stable category.

We  demonstrate that we can replicate some of these ideas in a more general setting. We need for this construction  a pair of additive categories and a triple of functors $(F,L,R)$ with $L$ and $R$ adjoint to $F$ satisfying some mild hypotheses. Sufficiently mild in fact that both of the best known (to an algebraist) triangulated categories, the homotopy category of an abelian category and the stable module category of a symmetric algebra, are examples of this construction. Of course, there is a formal definition of an algebraic triangulated category and it would be interesting to know what kinds of conditions on an algebraic triangulated category mean that it is (triangle equivalent to) one coming from a triple $(F,L,R)$. This is a slight generalization of the notion of a Frobenius pair.

In section \ref{functors} we explain how to use these functors to define a triangulated quotient. We then discuss the notion of transfer and give an equivalent definition of virtually $W$-projective that avoids  the transfer map (which will not, in our case, usually exist). We show that the virtually $W$-projective objects form a thick subcategory from this definition directly.  The idea then readily generalizes to our functor case. 

The central results of the paper are in section \ref{derived}, and subsequent sections.  We use $F$ to define  localizing subcategories of the homotopy categories. The localizations share many of the properties of the ordinary derived categories. The key ingredients  are the proof of the existence of $F$-projective resolutions of bounded above complexes. These satisfy certain extra conditions than being merely quasi-isomorphic. It should be noted that the resolutions of \cite{bok} are not sufficient for this purpose -- taking colimits does not in general preserve the properties we are interested in. We use our resolutions to identify the bounded relatively derived subcategories as the bounded homotopy categories of relatively projective objects. After we have identified some of the triangles in the relative derived category, we prove that the relatively stable category is a triangulated quotient of the bounded relative derived category.

\section{Preliminary results and definitions}

 We start by recalling some general nonsense results about categories. We  will use $\cata$ and $\catb$  to denote two exact categories. We refer the reader to \cite{keller} for the definition of an exact category. Suppose that  $F$ is an additive functor from $\cata$ to $\catb$.  We will show that if $F$ possesses adjoints, then under some assumptions these define a second  exact structure on $\cata$ with enough projectives/injectives.
 
\begin{lem}
 Suppose that $F$, $\cata$, and $\catb$ are as above. 
 \begin{enumerate}

 \item If $F$ has a left adjoint $L$ then $F$ is left exact. 

\item If $F$ has a right adjoint $R$ then $F$ is right exact.

\item If $F$ has both left and right adjoints and further if $\cata$ is abelian, then the counit of the adjunction  $\epsilon_X:LFX\to X$ is epic for all $X$  iff $FX=0$ implies $X=0$

\item If $F$ has both left and right adjoints and further if $\cata$ is abelian, then the unit of the adjunction  $\eta_X:X\to RFX$ is monic for all $X$  iff $FX=0$ implies $X=0$

\item Suppose that $F$ has a left  adjoint $L$, then the counit map $\epsilon_X:LFX\to X$ is epic for all $X$  iff $F$ is faithful.
\end{enumerate}
\end{lem}

The only statements not well known appear to be the $`FX=0$' results. The proof runs as follows for the first case.
\begin{proof} The `if' direction: consider the exact sequence $LFX \to X \to \mathrm{coker}(\epsilon_X)\to 0$ and apply the right exact $F$ (right exact by the existence of $R$) and we find $F\mathrm{coker}(\epsilon_X)=0$ hence the cokernel is zero as required. For the `only if 'direction, note that if $FX=0$ for some $X$, then there must be an epi from $LFX=0$ to $X$.\end{proof}
Remark: it would suffice for $F$ to be exact without  necessarily having both adjoints in the third and fourth statements; however we will only consider the case when $F$ is exact precisely because it has adjoints.

To see how we shall use these we recall a result of Happel \cite{happel} which we use to deduce that the constructions we make  yield triangulated quotients.  If $\cata$ is an exact category we say that an object $P$ is projective if the functor $(P,?)$ sends exact pairs to exact pairs. Dually an object $I$ is injective if $(?,I)$ sends exact pairs to exact pairs. We say $\cata$ has enough projectives  if for all $X$ there is an exact pair
\[ X' \to P \to X \]
with $P$ projective. The obvious condition characterizes enough injectives.
Our principal preliminary result is due to Happel \cite{happel}: 
\begin{thm} Suppose that $\mathcal{C}$ is an exact category with enough projective and injective objects. Further suppose that the class of projectives $\mathcal{P}$ coincides with the class of injectives then the quotient category
  \[ \underline{\mathcal{C}} := \frac{\mathcal{C}}{\mathcal{P}} \] is
  a triangulated category whose objects are those of $\mathcal{C}$ and
  $\mathrm{Hom}_{\underline{\mathcal{C}}}(X,Y)$ is the $\mathrm{Hom}$ set in $\mathcal{C}$
  modulo the relation $f\sim g$ iff $f-g$ factors through an object in
  $\mathcal{P}$.  The triangles in the quotient correspond to the exact pairs in $\mathcal{C}$.
  To define the shift functor, given $X$, pick some exact pair $X\to I \to X[1]$ with $I$ injective. This assignment is functorial on the quotient category; if we pick another exact pair and injective $X\to I' \to X'[1]$, then $X[1]$ is isomorphic to $X'[1]$ in the quotient.
\end{thm}
Let us give three key examples. For a longer account of instances of this we refer the reader to \cite{keller}.
\begin{exmp}[Exact Categories and Their Stable Categories]$\ $
\begin{itemize}
\item Let $\mathcal{C}$ be the category $\mathrm{mod}(\Lambda)$ of finitely generated modules for some finite dimensional symmetric algebra $\Lambda$. The exact structure is given by the set of all short exact sequences and the projective objects are precisely the projective modules. Since it is symmetric these are the same as injective modules. The triangulated category is the  stable module category.
\item If $H \leq G$ are finite groups, $k$ some algebraically closed field  then the subclass of short exact sequences of $kG$ modules that are split short exact upon restriction to $H$ endows $\mathrm{mod}(kG)$ with an exact structure. The projective objects are the direct summands of $kH$ modules induced up to $kG$.
\item Let $\mathcal{A}$ be some abelian category. Set $\mathrm{Ch}(\mathcal{A})$ to be the category of (co)chain complexes over $\mathcal{A}$. This is an additive category (even an abelian one but we do not need that). We define the exact structure using the short exact sequences of complexes that are split short exact sequences in each degree. The projective/injective objects are the contractible complexes, and the stable category is the homotopy category.
\end{itemize}\end{exmp}

\section{A triangulation of exact categories}\label{functors}
We proceed with the assumption that $\cata$ and $\catb$ are exact categories and that there is a functor $F:\cata\to \catb$ with left and right adjoints $L,R:\catb\to\cata$.  It is entirely possible that $\cata$ with its natural exact structure has  a triangulated quotient, however we will make absolutely no assumptions about whether $\cata$ comes with enough projectives, injectives, and even if it does whether they coincide. An example of an abelian category where the natural abelian structure does not give a triangulated quotient, but which has an exact substructure that does, is the category of chain complexes over an abelian category: triangles in the homotopy category do not correspond to short exact sequences of chain complexes.

We further suppose that all maps $\epsilon_X:LFX \to X$  and  $\eta_X:X\to RFX$ arising from the counit and unit of the adjunction are deflations and inflations respectively. If we were to think in terms of abelian categories instead then deflation should be replaced by epimorphism and inflation by monomorphism.  This is now sufficient to choose a subclass of exact pairs of $\cata$'s exact structure, and enough projective and injective objects with respect to this subclass.
\begin{defn}
Given functors and categories as above, we will define the $F$-projective objects to be the thick subcategory of all direct summands of the image of $L$
\[ \mathcal{P}:= \{ P \in \cata \ : \ P | LY\ \mathrm{for\ some}\ Y \in \catb\}.\]
Given some $X$ in $\cata$  \emph{the canonical $F$-projective cover} is the deflation  $\epsilon_X:LFX \to X$.  The translation functor is given by the exact pair
\[ \Omega_F(X)\to LFX \to X.\]
Dually we use $R$ and the unit $\eta$ to define $\mathcal{I}$ the class of $F$-injective objects, canonical $F$-injective hull and $\Omega_F^{-1}$ from the exact pair
\[ X\to FRX \to \Omega^{-1}(X).\]
Note that $\Omega$ and $\Omega^{-1}$ are only inverse functors in the quotient category.
Define $\Delta$ to be the  class of all exact  pairs in $\cata$ that  split upon applying $F$.
\end{defn}
\begin{prop} $\Delta$ is a second exact structure on $\cata$, and $\mathcal{P}$ and $\mathcal{I}$ are the classes of projective and injective objects with respect to $\Delta$. Further, there are enough projectives and injectives.
\end{prop}
\begin{proof}  Most of the axioms (see \cite{keller}) of an exact category are trivially satisfied. We explain only why $\Delta$ is closed under pull-backs.

To prove that we have (co)-cartesian squares we use the exactness of $F$. Given $p:Y\to Z$ an $F$-split deflation, and any morphism $W\to Z$ we can form the pull-back diagram since $\cata$ is exact and we obtain
\[ \xymatrix{ V\ar[d] \ar[r]^q & W \ar[d] \\ 
		Y \ar[r]^p&Z}\]
for some deflation $q$. If we apply $F$ then obtain another pull-back diagram for $Fp$ since $F$ is exact. Since $Fp$ is split, so is $Fq$. The dual statement for inflations holds.

Given $T$ in $\Delta$ and  $Y$ in $\catb$, it follows that  $(LY,T)$ = $(Y,FT)$. Since $FT$ is split $LY$ has the lifting property, hence is projective with respect to $\Delta$. Passing to summands we see $\mathcal{P}$ is the class of all projectives. Dually $\mathcal{I}$ is the class of injectives. That there are enough of each follows from the hypothesis on the unit and counit being inflations and deflations (they are obviously $F$-split).
\end{proof}
To summarize, we have the following theorem
\begin{thm}
Suppose that $\cata$ and $\catb$ are exact categories, $F$ is a functor from $\cata$ to $\catb$ with adjoints $L$ and $R$. If the $\epsilon_X:LFX \to X$ is always a deflation and $\eta_X:X\to RFX$ is always an inflation, and if further the subcategories of all summands of the images of $L$ and $R$ coincide, then by appeal to Happel's theorem there is a triangulated quotient category $\underline{A}_F$, and the triangles correspond to $F$-split exact pairs.
\end{thm}
And we also record the abelian version separately.
\begin{thm} If $\cata$ and $\catb$ are abelian categories and $F:\cata \to \catb$ is a faithful functor with left and right adjoints $L$ and $R$, then the class of short exact sequences that become split short exact upon applying $F$ is an exact structure on the underlying objects of $\cata$ with enough projectives and injectives. If the classes of injectives and projectives coincides then there is a triangulated quotient $\underline{A}_F$, and the triangles correspond to $F$-split short exact sequences.
\end{thm}
\begin{defn} 
When the hypothesis of the theorem are met we shall call the corresponding triangulated category the $F$-stable category. We will use $(?,?)_F$ for morphisms in the $F$-stable category.
\end{defn}
We end this section with two worked examples.
\begin{exmp}[Okuyama's Construction]\label{wstab}
Let $\cata = \catb = \mathrm{mod}(kG)$ for some finite group  $G$ and $k$ a field of characteristic dividing $|G|$ endowed with the exact structure of being an abelian  category. Fix $W$ some finite dimensional  module. With $F=W\otimes_k?$ we obtain Okuyama's relatively stable categories. Note this is trivial (in the sense that every object is projective) if $\mathrm{char}(k)$ does not divide $\mathrm{dim}(W)$.  In this special case we use $W$ rather than $W\otimes ?$ for $F$, and thus we talk of the $W$-stable category and $\Omega_W$ to agree with \cite{virtual}. If $W=kG$ we get the ordinary stable category. Note that as objects in the module category $\Omega_{kG}(X)$ and $\Omega(X)$, the usual Heller translate, are not isomorphic. $\Omega$ is only functorial at the level of the stable category, not the module category. $\Omega$ is normally defined using Krull-Schmidt, and we have avoided making any such assumptions about our exact categories.  
\end{exmp}
The second example is the homotopy category.
\begin{exmp}\label{homotopy}
Let $\mathcal{C}$ be an additive category, and $\cata$ the category of chain complexes over $\mathcal{C}$. In this case  we take $\catb$ to be the category of all $\mathbb{Z}$-graded objects over $\mathcal{C}$ so that a typical object is of the form 
\[ \coprod_{r \in \mathbb{Z}} X_r \]
We will suppose that the exact structure on $\cata$ is given by  the short sequences that are split exact in each degree. We take as $F$ the forgetful functor that  loses the differential information. We build the adjoint up piecewise. $L$ takes an object  $X_r$ in  degree $r$ and yields the complex that is zero everywhere except degrees $r$ and $r-1$ where it puts $X_r$ and the differential map is the identity. The right adjoint is similar except it uses degrees $r$ and $r+1$. The thick subcategory of summands of the images of $L$ and $R$ are the same: the contractible complexes. Note that we genuinely have an example where the adjoints $L$ and $R$ are different.
\end{exmp}

\section{Transfer, virtual projectivity  and the homotopy category}
A Frobenius pair  is a 2-tuple $(F,G)$ of functors with $F$ both left and right adjoint to $G$. This is of course a special case of our triple $(F,L,R)$. For a Frobenius pair, there is a well known concept of transfer, the prototype coming from the Frobenius pair of induction and restriction in representation theory.  If $H$ is a subgroup of $G$, then this map encodes much of the information about relative $H$-projectivity, in particular its image is the set of relatively projective maps, and $M$ is a relatively $H$-projective module if and only if the identity is in the image of $\mathrm{Tr}_H^G$ (the Higman Criterion). In \cite{virtual} Carlson, Peng and Wheeler examine  the relatively $W$-projective case, and prove similar results. That paper was one of the main inspirations for this work.  In their paper, the transfer map is defined as the composite
\[ \xymatrix{ \mathrm{Hom}_{kG}(M\otimes W, N\otimes W) \ar[r] &\mathrm{Hom}_{kG}(M,N\otimes W \otimes W^*) \ar[r] &\mathrm{Hom}_{kG}(M,N)}\]
where the second map is given by $(1_M\otimes \eta_N)_*$, where $\eta_N$ is the evaluation map arising from the unit of the adjunction.
The first of these maps is of course 
\[ (FX,FY) \cong (X,RFY)\]
and the second is
\[ (1\otimes \eta)_* : (X,LFY) \to (X,Y).\]
Since $L\cong R$  in the Frobenius case we can compose these morphisms.

We recall the definition from \cite{virtual}. 
\begin{defn} Suppose that $X$ is a $kG$-module, then $X$ is virtually $W$-projective if the induced map

\[ \mathrm{Tr}_W : \mathrm{Ext}_{kG}^n(X\otimes W,Y\otimes W) \to \mathrm{Ext}_{kG}^n(X,Y) \]
is surjective for all $n$ sufficiently large. 
\end{defn}
There  are, therefore, two obstructions. First we need a notion of transfer,  and secondly some Exts. 

The former of these can be sidestepped in special cases, such as the homotopy category where we can introduce a fudge factor, but not in general: suppose that $\cata$ is the category of chain complexes over some abelian category, and that $F$ is the forgetful functor to the category of $\mathbb{Z}$ graded objects with adjoints $L$ and $R$ as given in the last section. Recall that if $Y$ is a cochain complex then $RFY\cong Y\oplus Y[1]$.  It is a matter of elementary bookkeeping to show that the composition
\[ (X,Y) \to (FX,FY) \to (X,RFY) \cong (X,Y\oplus Y[1])\]
sends a map $g:X\to Y$ to $(g,dg)$. 
This leads us to the following notion: given $X$ and $Y$ cochain complexes, define the transfer map from $(FX,FY[-1])$ to $(X,Y)$ as the composite
\[ (FX,FY[-1]) \to (X,RFY[-1]) \cong (X,Y[-1]\oplus Y) \to (X,Y)\]
which takes degree zero {\em maps of graded objects} to chain maps. The image of  $\mathrm{Tr}(s)$ for $s$ a degree zero map of graded objects $FX\to FY[-1]$ is then $sd+ds$ which  is  an element  of the set of  maps homotopic to zero, and further, if $g$ is homotopic to zero, then $g=sd+ds$ for some $s$ a degree $1$ map of graded objects from $FX$ to $FY$ or equivalently a degree $0$ map of graded objects from $FX$ to $FY[-1]$, and it is in the image of the transfer map.

There is, however, an equivalent definition for virtually projective which does apply far more generally. Recall that if $kG$ is the group algebra of a finite group then   for $n \ge 1$ there is an isomorphism
\[ \mathrm{Ext}_{kG}^n(M,N)\cong (\Omega^n_{kG}(M),N)_{kG}.\] 
With the aid of the next lemma we can give an alternate definition for virtually $W$-projective modules  and show that they are  a thick triangulated subcategory of the $W$-stable category.
\begin{lem}
A module, $X$, is virtually $W$-projective if and only if  given a module $Y$, $(\Omega_{kG}^n(X),Y)_W=0$ for all $n$ sufficiently large. 
\end{lem}
\begin{proof}
The image of the transfer map in the module category is precisely the class of $W$-projective morphisms, see \cite{virtual} for example. If the transfer map is surjective from  $\mathrm{Ext}^n(W\otimes X, W\otimes Y)$ to Ext$^n(X,Y)$ then every morphism from $\Omega_{kG}(X) $ to $Y$ can be written (not necessarily uniquely) as $\alpha+\beta$ where $\alpha$ is $W$-projective and $\beta$ is projective (hence $W$-projective), and $(\Omega_{kG}^n(X),Y)_W = 0$. The converse is trivial.
\end{proof}
The key to an $\mathrm{Ext}$ free definition is the observation that  $\Omega_{kG}$ is actually a triangulated functor in the $W$-stable category. 
\begin{lem} The assignment $X\to \Omega_{kG}(X)$ is a triangulated endofunctor of the $W$-stable category.
\end{lem}
\begin{proof}
 In the module category $ \Omega_{kG}(X) \cong X\otimes \Omega_{kG}(k)$. Thus  $\Omega_{kG}(?)$ and $\Omega_{kG}(k)\otimes ?$ are  isomorphic as (additive)  functors  on the $W$-stable category. Tensoring over $k$ with a module is an exact functor on the module category, and it sends $W$-split short exact sequences to short exact sequences. Hence $\Omega_{kG}(?)$ is a triangulated functor on the $W$-stable category. 
\end{proof}
In general we have the following  proposition, which completes the $\mathrm{Ext}$ (and transfer) free definition of virtually $W$-projective modules, and proves they form a thick subcategory.
\begin{prop}
If $\mathcal{T}$ is a triangulated category, and $G$ a triangulated endofunctor, then 
\[ \mathcal{L} := \{ X \in \mathcal{T} :  \forall Y \in \mathcal{T}, (G^nX,Y)_\mathcal{T}=0 \textrm{ for all $n$ sufficiently large}\} \]
is a thick triangulated subcategory closed under direct sums 
\end{prop}
\begin{proof}
$\mathcal{L}$ is clearly closed under sums, summands and shifts because $G$ is a triangulated endofunctor. It only remains to argue that it is closed under triangles. If $T$ is a triangle in $\mathcal{T}$, and two terms are in $\mathcal{L}$, then for $n$ sufficiently large the  functor $(G^n ? ,Y)_\mathcal{T}$ applied to $T$ will give a long exact sequence and two out of every three consecutive terms vanish, hence they all vanish. 
\end{proof}
\begin{cor}
The class of virtually $W$-projective modules is a thick subcategory of the $W$-stable category.
\end{cor}

\begin{cor}\label{virt}
Suppose that $V$ is a $W$-projective module, then there is a (possibly trivial) thick subcategory $\mathcal{L}(V)$ of the $W$-stable category given by choosing $G=\Omega_V$, the Heller translate from the $V$-stable category.
\end{cor}
\begin{proof}
The proof that $\Omega_{kG}$ is a triangulated endofunctor of the $W$-stable category passes through with $V$ instead of $kG$.
\end{proof}

\noindent Remarks: \begin{enumerate}[(i)]
\item One can, if one felt so inclined, even choose $G$ to be the shift functor (or its inverse) of $\mathcal{T}$, or some combination of the two such as their sum. If $\mathcal{T}$ were the homotopy category of unbounded complexes, then the resulting virtually projective subcategories would be the bounded above, below, or totally subcategories.

\item It is worth spending a moment to consider the work that originally motivated this paper. We know that  in general it is too hard to classify all modules for a group algebra. However, we can and do borrow from topology and exploit the similarities between stable module categories and stable homotopy theory. Thus we are led to try to categorize modules up to some relation. And, owing to \cite{thick}, and \cite{balmer} we know that we should start by looking to exploit the extra tensor structure that modules have. A subcategory $\mathcal{I}$ is tensor closed if $X \in \mathcal{I}$ implies $X\otimes Y$ is in $\mathcal{I}$. Tensor closed thick  subcategories of the $\mathrm{stmod}(kG)$ are in one to one correspondence with specialization closed subsets of $V_G(k)$, the variety associated to the cohomology ring of $kG$. As of now, there is no corresponding object in the relatively stable case that plays the role of $V_G(k)$: relatively cohomology rings are not necessarily finitely generated. It is to be hoped that a framework will be developed that allows one to classify thick tensor ideals of relatively stable categories. This framework will have to allow for the fact that the thick subcategories of \ref{virt} are tensor ideal (simple exercise), and all trivial in the case of the ordinary stable module category. 
\end{enumerate}

\section{Basic structure of $D_F(\mathcal{A})$}\label{derived}
Suppose for the rest of the article that $\cata$ and $\catb$ are abelian. The derived category of $\cata$, when it exists, is the localization of K$(\cata)$ the homotopy category with the localizing class given by the acyclic complexes. We can construct a variation on the theme of a derived category that we label $D_F(\mathcal{A})$. We will  localize with respect to a multiplicative system. We will also prove the analogues  of  well known results that show bounded derived categories are equivalent to homotopy categories of bounded complexes of projectives. As is traditional when we talk about complexes we will mean cochain complexes: differentials raise degree.

We continue with the assumption that $F$ is faithful functor with left and right adjoints $L$ and $R$, and $\underline{\cata}_F$ is the the corresponding ($F$-stable category arising from the exact category $(\cata,\Delta)$. 
Since $F$ is exact it induces a functor from K$(\mathcal{A})$ the (unbounded homotopy category) of $\cata$ to K$(\catb)$.  Let $\mathcal{L}$ be the kernel of $F$ on K$(\cata)$, these are the acyclic complexes that split on applying $F$. This is certainly a localizing subcategory and multiplicative from an observation in \cite{koenig}[Proposition 2.5.1]  since it is also a null system as it is the kernel of a functor.
  \begin{defn} With $\cata$ and $\mathcal{L}$ as above define the localization
 \[  D_F(\cata):=\frac{\mathrm{K}(\cata)}{\mathcal{L}}\]
and the categories $D^-_F(\cata)$ is the localization of complexes in K whose terms are zero in all sufficiently high degrees (similarly for $+$ and $b$ instead of $-$). 
 \end{defn}
 Remark: we shall prove these bounded categories are full subcategories of $D_F(\cata)$ though they are not isomorphic  to objects with zero cohomology in sufficiently high (or low, or both)  degrees since it is possible for a complex to be acyclic yet not zero in $D_F(\cata)$. However these choices, allied to some careful bookkeeping allows us to prove most of the results we think of as natural.
In order to do anything though, we need to know when we may replace arbitrary complexes with complexes of $F$-projective and when we may not.  It is well known that there exists a quasi-isomorphic complex of relatively projective objects, \cite{hartshorne}, but that is not sufficient for our purposes: we need to show that there is not only a quasi-isomorphism, but one with an $F$-split cone.

\begin{defn} An $F$-projective resolution of a cochain complex $X$ is a complex whose terms are $F$-projective, with a map of complexes  \emph{in the category of complexes} $p:P \to X$ such that the cone of $p$ is zero in $D_F(\cata)$.
\end{defn}
In general the best we can do is to show that bounded above complexes have an $F$-projective resolution. Obviously in special cases we can  resolve all complexes, and if we were to impose some further categorical conditions on $\cata$ then we could perform some Bousfield localization argument. Without such restrictions we can prove:
\begin{thm} 
Let $X$ be a bounded above  complex in K$(A)$. Then $X$ has an $F$-projective resolution.
\end{thm}
\begin{proof}
We construct the complex $P$ explicitly. Let $X$ be a bounded above complex, then we may assume $X$ to be  zero in positive degrees.  Set $P^i$ to be zero for $i$ positive, and  $P^0$ be the canonical $F$-projective cover of $X^0$. We define,  inductively, the rest of the complex $P$: suppose we have defined 

\[ P^n\to P^{n+1} \to\ldots \to P^0\to 0\ldots\]
with maps from $P^r$ to $X^r$. Now consider the diagram
\[\xymatrix{
  & {\mathrm{ker}}(d^n_P) \ar[d] \\
 X^{n-1} \ar[r]^{d^{n-1}_X} & {\mathrm{ker}}(d^n_X) }\] 
from this, form the pull-back,
\[\xymatrix{
  Q \ar[r] \ar[d]& {\mathrm{ker}}(d^n_P) \ar[d] \\
 X^{n-1} \ar[r]^{d^{n-1}_X} & {\mathrm{ker}}(d^n_X) }\] 
and then define $P^n$ to be the canonical $F$-projective cover of $Q$:
\[\xymatrix{ P^n \ar[dr] \ar@{.>}[drr]^{d^{n-1}_P} \ar@{.>}[ddr]_{f_{n-1}}& & \\
 &  Q \ar[r] \ar[d]& {\mathrm{ker}}(d^n_P) \ar[d] \\
 &X^{n-1} \ar[r]^{d^{n-1}_X} & {\mathrm{ker}}(d^n_X) }\] 
We claim that $P \to X$ is an isomorphism on cohomology. If we momentarily grant ourselves this assumption, consider applying the functor Hom$_{\cata}(Y,?)$ for some $F$-projective object $Y$ in $\cata$ to obtain another map of complexes Hom$_{\cata}(Y,P) \to$ Hom$_{\cata}(Y,X)$. Since we chose everything carefully, this is also an isomorphism on cohomology since the functor  Hom$_{\cata}(Y,?)$ sends $F$-split short exact sequences to short exact sequences, and thus the complex Hom$_{\cata}(Y,\mathrm{cone}(p))$ is acyclic for any $F$-projective object $Y$ in $\cata$. This implies that Hom$_{\catb}(Z,F\mathrm{cone}(p))$ is acyclic for any $Z$ in $\catb$ and that cone$(p)$ is therefore split acyclic as required. Thus we need only prove it is an isomorphism on cohomology.
\begin{enumerate}[(i)]
\item $p$ is monic  (on cohomology): we know that Im$(d_P^{n-1}) \subset$ ker$(p)$, since p is a map of complexes. Let $[z] \in H^n(P)$ be realized by $z
  \in \mathrm{ker}(d^n_P)$, and suppose $p([z])=0$, then $p(z)$ is in
  Im$(d^{n-1}_X)$. Thus there is an element in the pull-back, $z'$, mapping to
  $z$ in $P^n$, and since the map from the canonical $F$-projective cover of the pull-back
  is an epimorphism $z$ lies in Im$(d^{n-1}_P)$, so ker$(p)\subset
  \mathrm{Im}(d^{n-1}_P)$, hence they are equal  and the map  
  \[ p: \frac{\mathrm{ker}(d^n_P)}{\mathrm{Im}(d^{n-1}_P)} \to
  \frac{\mathrm{ker}(d^n_X)}{\mathrm{Im}(d^{n-1}_X)}\]
is monic  as required.
\item $p$ is epic (on cohomology): let $[x]$ be in $H^n(X)$, then there is
  an element in the pull-back $(x,0)$ mapping to $x$ in $X^n$, and $0$ in
  $P^{n+1}$, and thus there is some $y \in P^n$ mapping to $x\in X^n$ and $0
  \in P^{n+1}$, ie $y \in \mathrm{ker}(d^n_P)$, so $p$ is an epimorphism. 
\end{enumerate}
\end{proof}
Naturally the dual statements about $F$-injective resolutions hold.

Let us recall some basic descriptions of maps in localizations. We may visualize maps in the derived category as roofs. A fraction $s^{-1}f:X
 \to Y$ is a roof:
\[\xymatrix{
   &  Z  &  \\
X \ar[ur]^f&                          & Y \ar[ul]^s}\] 
where $s$ and $f$ are maps in the homotopy category,and cone$(s)$ is $F$-split
acyclic. It is equivalent to the  fraction $t^{-1}g$ iff there is a commuting diagram
\[\xymatrix{
       &   Z\ar[d] &    \\
X\ar[ur]^f \ar[r] \ar[dr]^g & Z'   & Y \ar[ul]^s \ar[l]_u \ar[dl]^t\\
  & Z'' \ar[u] & }\]
with  the  mapping cone of the  map $u:Y \to Z'$ is $F$-split. A special case of this that we will need later is when $Y$ is a (bounded) complex of $F$-injective objects,  or $X$ a complex of $F$-projectives.

We will proceed by proving a series of straight forward results about maps in various $D_F(\cata)^\bullet$ that mimic those that we are used to in the usual case.
\begin{lem} Suppose that $P$ is a complex of relatively $F$-projective objects that is concentrated only in one degree, and that $S$ is an $F$-split acyclic complex, then 
\[ \mathrm{Hom}_{\mathrm{K}(\cata)}(S,P)=0 \]
\end{lem}
\begin{proof} It suffices to show this when $P$ is $LN$ for some object $N$ of $\catb$ thought of as a complex concentrated in degree zero. In this case
\[\mathrm{Hom}_{\mathrm{K}(\catb)}(S,LN) \cong \mathrm{Hom}_{\mathrm{K}(\cata)}(FS,N)=0\]
since $FS$ is a split exact sequence.
\end{proof}
\begin{lem}
Let $P$ be a bounded cochain complex and let $S$ be an $F$-split acyclic complex. If $P^i$ is  $F$-projective for all $i$, then
\[ \mathrm{Hom}_{\mathrm{K}(\cata)}(P,S)=0 \]
\end{lem}
\begin{proof} The class of objects for which this is true contains complexes of  $F$-projectives concentrated in one degree, and is closed under taking mapping cones, and hence contains all bounded complexes of  $F$-projectives.
\end{proof}
The next two proofs are easy exercises (c.f.~\cite{weibel}[Ch.~10]).
\begin{lem} Let $P$ be a bounded above complex of  $F$-projective objects and suppose that $f:P\to Y$ is a cochain map such that cone$(f)$ is an $F$-split acyclic complex, then $f$ is a split surjection.
\end{lem}
\begin{cor}[\cite{weibel}{\ Corollary 10.4.7}]
Let $P$ be a bounded above cochain complex of  $F$-projective objects, and $X$ any cochain complex, then
\[\mathrm{Hom}_{D_F(\cata)}(P,X)\cong \mathrm{Hom}_{\mathrm{K}(\cata)}(P,X)\]
\end{cor}
\begin{thm}
The category $D_F(\cata)^-$ is equivalent, as a triangulated category, to the homotopy category of bounded above complexes of $F$-projective objects.
\end{thm}
\begin{proof} Every complex  in $D_F(\cata)^-$  has an $F$-projective resolution . The homotopy category of complexes of bounded above $F$-projectives, call this $K^-(P(F))$, is a localizing  subcategory of $K^-(A)$. Thus, by general nonsense arguments \cite{weibel}{\ Section 10.3}, $D_F(\cata)^-$ is equivalent to the quotient category obtained by declaring the $F$-split acyclic objects of $K^-(P(F))$ to be zero. However, the only cochain complex in  this set of objects is the zero complex, hence  $D^-_F(A) \cong K^-(\mathcal{P})$.
\end{proof}
Remark: the dual statements about  $F$-injective objects have the corresponding dual proofs.

\section{Bounded categories are full subcategories}
It is well known that the usual categories $D^b(R), D^-(R)$, and $D^+(R)$ are full subcategories of $D(R)$ for $R$ a ring. This remains true for our $D^\bullet_F(\cata)$ in $D_F(\cata)$.  We shall need an initial result about how one may factor maps.
\begin{prop}
Let $M$ and $N$ be objects in the homotopy category. Suppose further that $N$ is bounded above and $f$ is a map from $M$ to $N$ and cone$(f)$ is an $F$ split acyclic. Then there is a bounded above complex $\hat{M}$ and  maps
\[g: \hat{M} \to M \]
\[ \hat{f}:\hat{M} \to N \]
such that cone$(\hat{f})$ is an $F$ split acyclic complex and $\hat{f}=fg$.
\end{prop}
\begin{proof}
Without loss of generality we may suppose that $N$ lies only in negative degree so that $N^i\cong 0$ for all $i>0$.  Let $M'$ be the `sensible' truncation of $M$ in degree 1, that is $\hat{M}^i= M^i$ for $i \leq 1$, $\hat{M}^2=\mathrm{Im}(d^1_M)$ the image of the differential in degree $0$ of $M$, and $\hat{M}^i \cong 0$ for $i >2$, and then $g$ is the inclusion of $\hat{M}$ into $M$ and $\hat{f}$ is the composite $gf$. If we examine cone$(f)$ then it becomes clear that these objects satisfy the conditions of the theorem:
\[ \mathrm{cone}(f) = \ldots N^{-1} \oplus M^0 \to \underbrace{N^0\oplus M^1}_{\mathrm{degree}\ 0} \to M^2 \to M^3 \ldots \]
which we know to be $F$ split acyclic, and in particular we can decompose $F(\mathrm{cone}(f))$ into the direct sum of two split acyclics
\[ \ldots FN^{-1} \oplus FM^0 \to FN^0\oplus FM^1 \to \mathrm{Im}(Fd^1_M) \to 0 \to 0 \ldots \]
and 
\[  \ldots 0  \to  \mathrm{ker}(Fd_2) \to FM^3 \ldots \]
the former of these is exactly cone$(\hat{f})$ and the second shows cone$(g)$ to be $F$ split exact  and we see all the requirements of the theorem are met. 
\end{proof}
This immediately gives
\begin{cor}
The inclusions of $D^-_F(\cata)$ in $D_F(\cata)$,  $D^b_F(\cata)$ in $D^+_F(\cata)$, $D^-_F(\cata)$ into $D_F(\cata)$ and $D^b_F(\cata)$ into $D^-_F(\cata)$ are all fully faithful. Moreover $D^-_F(\cata)$ is equivalent to the full subcategory of complexes which are $F$-split acyclic in sufficiently high degrees, and the analogous statements hold for $D^b_F(\cata)$ and $D^+_F(\cata)$.
\end{cor}
\begin{proof} We prove only the first of these. Let $N$ and $N'$ be two bounded above complexes in $D_F(\cata)$. Then a morphism is a roof
\[\xymatrix{ & M \ar[dr] \ar[dl]^f& \\
           N&&N' }\]
with cone$(f)$ an $F$ split acyclic. If we keep the  notation of the last proposition then we can complete the diagram
\[\xymatrix{  & M \ar[dl]_f \ar[drrr]& & {\hat{M}} \ar[dlll]^{\hat{f}} \ar[dr] \\
             N &&& &N' }\]
to 
\[\xymatrix{ && {\hat{M}} \ar[dl]^g \ar[dr]^1 && \\
 & M \ar[dl]_f \ar[drrr] & &{\hat{M}} \ar[dlll]^{\hat{f}} \ar[dr] \\
             N &&& &N' }\]
and in particular every map between bounded above complexes is equivalent to a map coming from the inclusion functor.

Essentially the same argument shows that the inclusion functor is faithful too.
For if $hf^{-1}: X \to Y$ is a map in $D^-_F(\cata)$ that becomes zero in $D_F(\cata)$ then we are saying that a the diagram
\[\xymatrix{  & Z \ar[dl]_f \ar[drrr]_h & & Z \ar[dlll]^f \ar[dr]^0 \\
             X &&& &Y }\]
where all objects are bounded above can be completed to a commuting diagram
\[\xymatrix{ && Z' \ar[dr] \ar[dl] && \\
               & Z \ar[dl]_f \ar[drrr]_h & & Z \ar[dlll]^f \ar[dr]^0 \\
             X &&& &Y }\]
where $Z'$ is a possibly unbounded complex and the maps from $Z'$ are invertible in $D_F(\cata)$. But we know that we can replace $Z'$ by a complex  bounded above from argument we gave for  the fullness of the inclusion, hence the map was zero in $D^-_F(\cata)$ already.

 It should be clear that $D^-_F(\cata)$, say,  is the full subcategory objects $X$ that fit into a distinguished triangle $X^- \to X^+ \to X$ with $X^-$ a bounded above complex and $X^+$ an $F$-split acyclic.   
\end{proof}

\section{Triangles in $D_F(\cata)$}
In order to work with  $D^{\bullet}_F(\cata)$, and show that $D^b_F(\cata)$ has the stable category as a quotient,  we should discuss the triangulated structure of $D_F(\cata)$, or more precisely we should describe the triangles. We know that in the usual derived category the distinguished triangles arise from short exact sequences of cochain complexes, that is to say, if
\[0 \to X \to Y \to Z \to 0 \]
is a short exact sequence of complexes, then there is a triangle in the derived category (but not necessarily the homotopy category)
\[ X \to Y \to Z \to X[1] \]
The natural generalization of this is
\begin{prop}\label{triangles} If there is a short exact sequence of complexes
\[ X \to Y \to Z\]
 such that 
\[ FX\to FY \to FZ \]
is a split short exact sequence (recall F is exact), then there  is a distinguished triangle 
\[X \to Y \to Z \to X[1] \]
in $D_F(\cata)$.
\end{prop}
\begin{proof} We know that there is a diagram of maps of complexes
\[\xymatrix{  X \ar[r]^f & Y\ar[r]^g & Z \\
              X \ar[u]^1 \ar[r] & Y \ar[u]^1 \ar[r]& {\mathrm{cone}}(f) \ar[u]^{\phi}}\]
and that $\phi$ is a quasi-isomorphism. We need only show that cone$(\phi)$ is $F$-split, as it then becomes an isomorphism in $D_F(\cata)$.
\[ F\mathrm{cone}(\phi)^n = FX^{n+2}\oplus FY^{n+1} \oplus FZ^n \]
and the differential (abusing $d$ as usual) is
\[ d = \left( \begin{array}{ccc} Fd & 0 &0\\
                            Ff & -Fd & 0\\
                           0& Fg  & Fd \end{array} \right) \]
But $Ff$ and $Fg$ are split maps, let the splitting maps be $f'$ and $g'$ respectively then
\[ s = \left( \begin{array}{ccc}  0 & f' &0 \\
                             0&  0&  g'\\
                             0&0& 0\end{array} \right)\]
satisfies
\[sd+ds= \left( \begin{array}{ccc} f'Ff  & -f'd +df' & 0\\
                                   0 & g'Fg+Fff' &0\\
                                   0 & 0 & Fgg' \end{array} \right) \]
that is to say
\[sd+ds= \left( \begin{array}{ccc} 1_{X^{n+2}}& 0& 0\\
                                   0 & 1_{Y^{n+1}} &0\\
                                   0 & 0 & 1_{Z^n}\end{array} \right) \]
as $f'$ is a chain map (commutes with differentials).
\end{proof}

\section{The $F$-stable category as a quotient}
 We are now in a position to prove a theorem which is akin to a result originally due to Buchweitz in his unpublished manuscript, \cite{buch}. The proof, given the existence of $F$-projective resolutions, is essentially the same as that given by Rickard in \cite{stable}. We include it for completeness.
\begin{thm}
Let $D^b_F$ be the bounded relative derived category considered as the full subcategory of complexes in $D_F(\cata)$ that are $F$-split acyclic in sufficiently high and low degrees. If $K(\mathcal{P})$ denotes the full thick triangulated subcategory of complexes whose terms are $F$-projective modules lying in finitely many degrees, then the inclusion of $\mathcal{A}$  into complexes concentrated in degree zero induces an equivalence
\[ S:\underline{\cata}_F \to \frac{D^b_F(\cata)}{K(\mathcal{P})} \]
\end{thm}
\begin{proof} The natural inclusion of $\cata$ sends $F$-projective objects to zero in the quotient, hence factors through the  stable category, inducing $S$. A distinguished triangle in the  stable category is equivalent to an $F$-split short exact sequence, which gives under the inclusion, an $F$-split short exact sequence of cochain complexes, and hence  a triangle in $D_F(\cata)$, so $S$ is a triangulated functor.

We next show $S$ is surjective on (isomorphism classes of) objects. Let $X$ be an element of $D^b_F(\cata)$. We may suppose $X$ is only non-zero in finitely many degrees, and then we may find an $F$-projective resolution $P$ of $X$. We may assume that $P$ is zero in positive degrees, and that it is $F$-split acyclic in degrees less than $r$ for some $r$.  Let $P'$ be the naive truncation of $P$ to the terms in degrees less than $r$. $P$ and $P'$ are isomorphic in the quotient $D^b_F(\cata)/K(\mathcal{P})$ since the cone of the inclusion is a bounded complex of $F$-projectives. We may extend $P'$  to $P''$ by taking an $F$-injective resolution of the cohomology of $P'$ lying in degree $r$, and terminate it in degree zero. The map from $P''$ to $P'$ is an isomorphism in the derived category, and $P''$ is the $F$-projective resolution of some module in degree 0, hence isomorphic to an object in the image of $S$.

We must show $S$  is full. A map $X\to Z$ in the quotient category $D^b_F(\cata)/K(\mathcal{P})$ is a roof of morphisms in $D^b_F(\cata)$:
\[\xymatrix{ & Y & \\
            X \ar[ur]&& Z \ar[ul] }\]
            which fits into a diagram
\[\xymatrix{  P & & \\        
                &Y\ar[ul]&\\
             X\ar[ur] && Z\ar[ul]} \]
with arrows morphisms in $D^b_F(\cata)$ and $P$ a bounded complex of $F$-projectives. If $X$ and $Z$ are complexes concentrated in degree zero we will show that the map $X \to P$  is zero in $D_F(\cata)$ and thus that there is a map in $D^b_F(\cata)$ from $X$ to $Z$.  As before we may take $Y$ to be a complex of  $F$-projective objects in negative degrees, and hence $P$ can be taken to lie in strictly negative degree. Thus there are no maps from $X$ to $P$ and there is a map in homotopy category between the  $F$-projective resolutions of $X$ and $Y$. Clearly, then, $S$ is full.  Equally, it follows that a complex concentrated in degree zero is isomorphic to an object of $K(\mathcal{P})$ if and only if the degree zero term  is  $F$-projective. 

By a general nonsense argument about triangulated functors, this shows $S$ is faithful:  suppose that $Sf=0$, and consider the distinguished triangle
\[\xymatrix{ X \ar[r]^f& Y \ar[r]^g & Z \ar[r] & X[1] }\]
 Applying $S$ to this triangle we obtain (one isomorphic to)
\[\xymatrix{  SX \ar[r]^0  & SY \ar[r] &
  SY \oplus SX[1] \ar[r] & SX[1]  }\] 
since $Sf=0$. Let $h'$ be the splitting map 
\[ SY \oplus SX[1] \to SY \]
such that  $h'Sg = \mathrm{Id}_{SY}$. Since $S$ is full, there is an $h$
such that $h'=Sh$, i.e.~$Shg = \mathrm{Id}_{SY}$. Consider the triangle
\[\xymatrix{ Y \ar[r]^{hg}& Y \ar[r]& V \ar[r] & Y[1] }\]
 If we apply $S$ to this, then $SV$ must vanish. Hence $V$ is
isomorphic to zero, and $hg$ is an isomorphism. We conclude that  $g$ is a split monomorphism, and $f$ must be zero, as we were required to show.
\end{proof}

\section{Acknowledgements}
This work was carried out whilst a PhD student at the University of Bristol under the kind and knowledgeable supervision of Jeremy Rickard to whom I am greatly indebted. The position was funded by the EPSRC. I also benefited from the advice and encouragement of Joe Chuang and Thorsten Holm. I would also like to thank Mike Wemyss for pointing out a simplification of the omnibus lemma for $F$ and its adjoints in the case of abelian categories. I would also like to thank the referee for suggesting improvements to the manuscript and pointing out that the original hypotheses of example \ref{homotopy} were overly restrictive.

\end{document}